\documentclass[10pt]{article}
\usepackage{setspace}

\usepackage{amsthm,amsmath,amssymb}

\newtheorem{prp}{Proposition}
\newtheorem{lmm}{Lemma}
\newtheorem{thr}{Theorem}
\newtheorem{example}{Example}

\begin{document}

\begin{center}
{\Large \bf Cycle structures of autotopisms of the Latin squares
of order up to 11.}

\vspace{0.5cm}

{\large Falc\'on, R. M.}

{\small Department of Geometry and Topology. \\
Faculty of Mathematics. University of Seville. \\
41080 - Seville (Spain). \\
E-mail: {\em rafalgan@us.es}}
\end{center}

\begin{center} {\large \bf Abstract} \end{center}

The cycle structure of a Latin square autotopism
$\Theta=(\alpha,\beta,\gamma)$ is the triple
$(\mathbf{l}_{\alpha},\mathbf{l}_{\beta},\mathbf{l}_{\gamma})$,
where $\mathbf{l}_{\delta}$ is the cycle structure of $\delta$,
for all $\delta\in\{\alpha,\beta,\gamma\}$. In this paper we study
some properties of these cycle structures and, as a consequence,
we give a classification of all autotopisms of the Latin squares
of order up to 11.

\vspace{0.2cm}



\vspace{0.3cm}

\section{Introduction}

A {\em quasigroup} {\bf \cite{Albert43}} is a nonempty set $G$
endowed with a product $\cdot$, such that if any two of the three
symbols $a,b,c$ in the equation $a\cdot b=c$ are given as elements
of $G$, the third is uniquely determined as an element of $G$. It
is equivalent to say that $G$ is endowed with left and right
division. Two quasigroups $(G,\cdot)$ and $(H,\circ)$ are {\em
isotopic} {\bf \cite{Bruck44}} if there are three bijections
$\alpha,\beta,\gamma$ from $H$ to $G$, such that $\gamma (a\circ
b) = \alpha(a)\cdot \beta(b)$, for all $a,b\in H.$ The triple
$\Theta=(\alpha,\beta,\gamma)$ is called an {\em isotopism} from
$(H,\circ)$ to $(G,\cdot)$. The multiplication table of a
quasigroup is a Latin square. A {\em Latin square} $L$ of order
$n$ is a $n \times n$ array with elements chosen from a set
$N=\{x_1,...,x_n\}$, such that each symbol occurs precisely once
in each row and each column. The set of Latin squares of order $n$
is denoted by $LS(n)$. The calculus of the number of Latin squares
of order $n$ is an open problem. However, this number is known up
to order $11$ {\bf \cite{McKay05a}}. A general overview of Latin
squares and their applications can be seen in {\bf \cite{Denes91}}
or {\bf \cite{Laywine98}}.

\vspace{0.15cm}

Throughout this paper, we will consider $N=\{0,1,...,n-1\}$ and
$S_n$ will denote the symmetric group on $N$. The {\em cycle
structure of a permutation} $\delta\in S_n$ is the sequence
$(\mathbf{l}_1,\mathbf{l}_2,...,\mathbf{l}_n)$, where
$\mathbf{l}_i$ is the number of cycles of length $i$ in $\delta$.
For a given $\delta\in S_n$, define the set of its fixed points by
$Fix(\delta)=\{i\in N : \delta(i)=i\}$. If
$L=\left(l_{i,j}\right)\in LS(n)$, the {\em orthogonal array
representation of} $L$ is the set of $n^2$ triples
$\{(i,j,l_{i,j}): i,j\in N\}$. The previous set is identified with
$L$ and so, it is written $(i,j,l_{i,j})\in L$, for all $i,j\in
N$. Moreover, since $L$ is the multiplication table of a
quasigroup, then distinct triples of $L$ never agree in more than
one element.

\vspace{0.15cm}

An {\em isotopism} of a Latin square $L\in LS(n)$ is a triple
$\Theta=(\alpha,\beta,\gamma)\in \mathcal{I}_n=S_n\times S_n\times
S_n$. So, $\alpha,\beta$ and $\gamma$ are permutations of rows,
columns and symbols of $L$, respectively. The resulting square
$L^{\Theta}$ is also a Latin square and it is said to be {\em
isotopic} to $L$. In particular, if $L=\left(l_{i,j}\right)$, then
$L^{\Theta}=\{(i,j,\gamma^{-1}\left(l_{\alpha(i),\beta(j)}\right):
i,j\in N\}$. If $L_1$ and $L_2$ are two distinct Latin squares of
order $n$, then $L_1^{\Theta}\neq L_2^{\Theta}$. If
$\alpha=\beta=\gamma$, the isotopism is an {\em isomorphism}. If
$\gamma=\epsilon$, the identity map on $N$, $\Theta$ is called a
{\em principal isotopism}. An isotopism which maps $L$ to itself
is an {\em autotopism}. Moreover, if it is an isomorphism, then it
is called an {\em automorphism}. If its permutations are $n$
cycles, then $L$ is said to be {\em diagonally cyclic}.  Indeed,
diagonally cyclic Latin squares of even order do not exist {\bf
\cite{Wanless04}}. $(\epsilon,\epsilon,\epsilon)$ is called the
{\em trivial autotopism}. The stabilizer subgroup of $L$ in
$\mathcal{I}_n$ is its {\em autotopism group},
$\mathcal{U}(L)=\{\Theta\in \mathcal{I}_n:L^{\Theta}=L\}$. For a
given $L\in LS(n)$, $\Theta=(\alpha,\beta,\gamma)\in
\mathcal{U}(L)$ and $\sigma\in S_3$, it is verified that
$(\pi_{\sigma(0)}(\Theta),\pi_{\sigma(1)}(\Theta),\pi_{\sigma(2)}(\Theta))\in
\mathcal{U}(L^{\sigma})$, where $\pi_i$ gives the $(i+1)^{th}$
component of $\Theta$, for all $i\in\{0,1,2\}$. For a given
$\Theta\in\mathcal{I}_n$, the set of all Latin squares $L$ such
that $\Theta\in \mathcal{U}(L)$ is denoted by $LS(\Theta)$. The
cardinality of $LS(\Theta)$ is denoted by $\Delta(\Theta)$.
Specifically, the computation of $\Delta(\Theta)$ for any
isotopism $\Theta\in \mathcal{I}_n$ is at the moment an open
problem having relevance in secret sharing schemes related to
Latin squares {\bf \cite{LongCycles}}.

\begin{figure}[htb]\label{Fig1}{\small
\begin{center} $\begin{cases} L_1=\left(\text{\tiny \begin{tabular}{cccc}
  0 & 1 & 2 & 3 \\
  1 & 0 & 3 & 2 \\
  2 & 3 & 0 & 1 \\
  3 & 2 & 1 & 0
\end{tabular}}\right)  \\ \ \\ \Theta=((0\ 1)(2\ 3), (1\ 2),\epsilon)\end{cases}
\Rightarrow \ L_1^{\Theta}=\left(\text{\tiny \begin{tabular}{cccc}
  1 & 3 & 0 & 2 \\
  0 & 2 & 1 & 3 \\
  3 & 1 & 2 & 0 \\
  2 & 0 & 3 & 1
\end{tabular}}\right)$
\caption{Isotopism permuting $1^{st}$ with $2^{nd}$ and $3^{rd}$
with $4^{th}$ rows and $2^{nd}$ with $3^{rd}$ columns.}
\end{center}}
\end{figure}

\vspace{0.15cm}

The following result gives some necessary conditions of the
possible non-trivial Latin square autotopisms:

\begin{thr}[McKay, Meynert and Myrvold {\bf \cite{McKay05}}]
\label{thr1} Let $L\in LS(n)$. Every non-trivial
$\Theta=(\alpha,\beta,\gamma)\in\mathcal{U}(L)$ verifies one of
the following assertions:
\begin{enumerate}
\item[a)] $\alpha,\beta,\gamma$ have the same cycle structure with
at least one and at most $\lfloor \frac n2 \rfloor$ fixed points.
\item[b)] One of $\alpha,\beta,\gamma$ has at least one fixed
point and the other two have the same cycle structure without
fixed points.
\item[c)] None of $\alpha,\beta,\gamma$ has fixed points. \hfill $\Box$ \vspace{0.5cm}
\end{enumerate}
\end{thr}

The classification given in the previous theorem depends on the
cycle structures of the permutations of each Latin square
autotopism and on their fixed points. In this paper, we are
interested in giving a complete catalogue with all the possible
cycles structures of any autotopism of a Latin square of order up
to 11. This catalogue seems to be useful to study the open problem
of the calculus of the number $\Delta(\Theta)$. Specifically, we
prove in Section 3 that the number of Latin squares having a given
isotopism $\Theta\in\mathcal{I}_n$ in its autotopism group only
depends on the cycle structure of $\Theta$.

\vspace{0.25cm}

The structure of the paper is the following: in Section 2, some
general results about Latin square autotopisms are reviewed. In
Section 3, we define the cycle structure of a Latin square
autotopism and we study several of its properties. All these
properties have been implemented in a computer program to give in
Section 4 the classification of all autotopisms of the Latin
squares of order up to 11.

\section{Some general results}

Every permutation of $S_n$ can be written as the composition of
pairwise disjoint cycles. So, from now on, for a given
$\Theta=(\alpha,\beta,\gamma)\in\mathcal{I}_n$, we will consider
that, for all $\delta\in\{\alpha,\beta,\gamma\}$:

\vspace{0.15cm}

$\hspace{3.25cm} \delta=C^{\delta}_0\circ C^{\delta}_1\circ ...
\circ C^{\delta}_{k_{\delta}-1}, \hspace{3.25cm} (1)$

\noindent where:
\begin{enumerate}
\item[i)] For all $i\in \{0,1,...,k_{\delta}-1\}$, one has
$C^{\delta}_i=\left(c_{i,0}^{\delta}\ c_{i,1}^{\delta}\ ...\
c_{i,\ \lambda_i^{\delta}-1}^{\delta}\right)$, with
$\lambda_i^{\delta}\leq n$ and $c_{i,0}^{\delta}=\min_j
\{c_{i,j}^{\delta}\}$.
\item[ii)] $\sum_i \lambda_i^{\delta}=n$.
\item[iii)] For all $i,j\in \{0,1,...,k_{\delta}-1\}$, one has
$\lambda_i^{\delta}\geq \lambda_j^{\delta}$, whenever $i\leq j$.
\item[iv)] Given $i,j\in \{0,1,...,k_{\delta}-1\}$, with $i<j$ and
$\lambda_i^{\delta}=\lambda_j^{\delta}$, one has
$c_{i,0}^{\delta}<c_{j,0}^{\delta}$.
\end{enumerate}

Specifically, the following result is verified:

\begin{prp}
\label{prp1} Let $\Theta=(\alpha,\beta,\gamma)\in \mathcal{I}_n$
be a non-trivial isotopism. If one of the permutations
$\alpha,\beta$ or $\gamma$ is equal to $\epsilon$, then
$\Delta(\Theta)$ $>0$ only if the other two permutations have the
same cycle structure with all their cycles of the same length and
without fixed points.
\end{prp}

{\bf Proof.} Let $\Theta=(\alpha,\beta,\gamma)\in \mathcal{I}_n$
be such that $\Delta(\Theta)>0$ and let us consider
$L=(l_{i,j})\in LS(\Theta)$. If one of the permutations
$\alpha,\beta$ or $\gamma$ is equal to $\epsilon$, then we are in
case (b) of Theorem \ref{thr1} and, therefore, the other two
permutations must have the same cycle structure without fixed
points. Now, we must prove that all the cycles of these two
permutations have the same length. To do it, since rows, columns
and symbols have an interchangeable role in the study of Latin
squares, it is enough to study the case $\alpha=\epsilon$, being
equivalent the proof when $\beta=\epsilon$ or $\gamma=\epsilon$.
Thus, $\beta$ and $\gamma$ have the same cycle structure without
fixed points. Specifically, $k_{\beta}=k_{\gamma}$. Let us suppose
that there exist $r,s\in\{0,1,...,k_{\beta}-1\}$ such that
$\lambda_r^{\beta}\neq\lambda_s^{\gamma}$. Now, let $a\in N$ be
such that $l_{a,\ c_{r,0}^{\beta}}=c_{s,0}^{\gamma}$. If
$\lambda_r^{\beta}>\lambda_s^{\gamma}$, then:
$$l_{a,\ c_{r,0}^{\beta}}=c_{s,0}^{\gamma}=c_{s,\lambda_s^{\gamma}\ (mod\ \lambda_s^{\gamma})}^{\gamma}=
l_{a,\ c_{r,\lambda_s^{\gamma}}^{\beta}},$$ which is a
contradiction with being $L$ a Latin square. Otherwise, if
$\lambda_r^{\beta}<\lambda_s^{\gamma}$, then:
$$c_{s,0}^{\gamma}=l_{a,\ c_{r,0}^{\beta}}=
l_{a,\ c_{s,\lambda_r^{\beta}\ (mod\
\lambda_r^{\beta})}^{\beta}}=c_{s,\lambda_r^{\beta}}^{\gamma},$$
which is a contradiction with the conditions $(1)$ imposed at the
beginning of this section. Therefore, it must be that
$\lambda_r^{\beta}=\lambda_s^{\gamma}$, for all
$r,s\in\{0,1,...,k_{\beta}-1\}$. \hfill $\Box$ \vspace{0.5cm}

From now on, for a given $\delta\in\{\alpha,\beta,\gamma\}$ and
$i\in \{0,1,...,k_{\delta}-1\}$, we will write $a\in C^{\delta}_i$
if there exists $j\in\{0,1,...,\lambda_i^{\delta}-1\}$ such that
$a=c_{i,j}^{\delta}$. The following result is verified:

\begin{thr}\label{thr2} Let $L=(l_{i,j})\in LS(n)$ and
$\Theta=(\alpha,\beta,\gamma)\in\mathcal{U}(L)$ and let us
consider $r\in\{0,1,...,k_{\alpha}-1\}$ and
$s\in\{0,1,...,k_{\beta}-1\}$. Let us denote
$m=l.c.m.(\lambda_r^{\alpha},\lambda_s^{\beta})$. Now, for a given
$a\in C_r^{\alpha}$ and $b\in C_s^{\beta}$, let
$t\in\{0,1,...,k_{\gamma}-1\}$ be such that $l_{a,b}\in
C_t^{\gamma}$. Then, it is verified that:
\begin{enumerate}
\item[i)] $\lambda_t^{\gamma}$ divides $m$.
\item[ii)] $\lambda_t^{\gamma}$ does not divide any multiple of
$\lambda_r^{\alpha}$ smaller than $m$.
\item[iii)] $\lambda_t^{\gamma}$ does not divide any multiple of
$\lambda_s^{\beta}$ smaller than $m$.
\item[iv)] If $g.c.d.(\lambda_r^{\alpha},\lambda_s^{\beta})=1$, then $\lambda_t^{\gamma}=m$.
\end{enumerate}
\end{thr}

{\bf Proof.} Let $u\in \{0,1,...,\lambda_r^{\alpha}-1\}$, $v\in
\{0,1,...,\lambda_s^{\beta}-1\}$ and $w\in
\{0,1,...,\lambda_t^{\gamma}-1\}$ be such that
$a=c_{r,u}^{\alpha}$, $b=c_{s,v}^{\beta}$ and
$l_{a,b}=c_{t,w}^{\gamma}$, respectively. Since $\Theta\in
\mathcal{U}(L)$, we obtain that $\lambda_t^{\gamma}$ divides $m$,
because it must be that:
$$c_{t,w}^{\gamma}=l_{a,b}=l_{c_{r,u}^{\alpha},\
c_{s,v}^{\beta}}=l_{c_{r,u + m\ (mod \
\lambda_r^{\alpha})}^{\alpha},\ c_{s,v+m \ (mod\
\lambda_s^{\beta})}^{\beta}}=c_{t,w+m\ (mod\
\lambda_t^{\gamma})}^{\gamma}.$$

Now, let us suppose that $\lambda_r^{\alpha}\neq
\lambda_s^{\beta}$. Then, we see that $\lambda_t^{\gamma}$ does
not divide any multiple $h$ of $\lambda_r^{\alpha}$ smaller than
$m$:
$$c_{t,w}^{\gamma}=l_{c_{r,u}^{\alpha},\
c_{s,v}^{\beta}}=l_{c_{r,u+ h\ (mod \
\lambda_r^{\alpha})}^{\alpha},\ c_{s,v}^{\beta}}\neq l_{c_{r,u+ h\
(mod \ \lambda_r^{\alpha})}^{\alpha},\ c_{s,v+h \ (mod\
\lambda_s^{\beta})}^{\beta}}=$$ $$=c_{t,w+h\ (mod\
\lambda_t^{\gamma})}^{\gamma}.$$

In a similar way, it can be obtained that $\lambda_t^{\gamma}$
does not divide any multiple of $\lambda_s^{\beta}$ smaller than
$m$.

Finally, if $g.c.d.(\lambda_r^{\alpha},\lambda_s^{\beta})=1$, then
$m=\lambda_r^{\alpha}\cdot \lambda_s^{\beta}$. Let us suppose that
$\lambda_t^{\gamma}<m$. By keeping in mind assertions (ii) and
(iii), since $g.c.d.(\lambda_r^{\alpha},\lambda_s^{\beta})=1$,
there must exist two distinct primes $p,q\in [m]$ such that $p$
divides $\lambda_r^{\alpha}$, $q$ divides $\lambda_s^{\beta}$ and
$\lambda_t^{\gamma}$ divides $\frac m{p\cdot q}$. Specifically,
$\lambda_t^{\gamma}$ divides $\frac mp$, which is a multiple of
$\lambda_s^{\beta}$. It is a contradiction with assertion (iii)
and, therefore, it must be that  $\lambda_t^{\gamma}=m$. \hfill
$\Box$ \vspace{0.5cm}

\section{\hspace{-0.5cm} Cycle structures of Latin square autotopisms}

From now on, for a given $n\in\mathbb{N}$, we will denote the set
$\{1,2,...,n\}$ by $[n]$. So, let $\Theta=(\alpha,\beta,\gamma)\in
\mathcal{I}_n$ and let us define, for all
$\delta\in\{\alpha,\beta,\gamma\}$ and $r\in [n]$:
$$\mathbf{l}_r^{\delta}=\sharp\{i\in\{0,1,...,k_{\delta}-1\}:\lambda^{\delta}_i=r\},$$
where $\sharp$ denotes the cardinality of the corresponding set.
Then, let us consider, for all $\delta\in\{\alpha,\beta,\gamma\}$:
$$\mathbf{l}_{\delta}=(\mathbf{l}^{\delta}_1,\mathbf{l}^{\delta}_2,...,\mathbf{l}^{\delta}_n).$$

The triple
$(\mathbf{l}_{\alpha},\mathbf{l}_{\beta},\mathbf{l}_{\gamma})$
will be called the {\em cycle structure of} $\Theta$. The set of
all autotopisms of the Latin squares of order $n$ having the cycle
structure
$(\mathbf{l}_{\alpha},\mathbf{l}_{\beta},\mathbf{l}_{\gamma})$
will be denoted by $\mathcal{I}_{n}(\mathbf{l}_{\alpha},\
\mathbf{l}_{\beta},\ \mathbf{l}_{\gamma})$.

\vspace{0.25cm}

Some immediate properties of the cycle structure of an isotopism
are given in the following:

\begin{lmm} \label{lmm1} Let $\Theta\in \mathcal{I}_{n}(\mathbf{l}_{\alpha},\ \mathbf{l}_{\beta},\
\mathbf{l}_{\gamma})$. Then, for all
$\delta\in\{\alpha,\beta,\gamma\}$, it must be that:
\begin{enumerate}
\item[a)] $\sum_{r\in [n]} \mathbf{l}_r^{\delta}=k_{\delta}$.
\item[b)] $\sum_{r\in [n]} r\cdot \mathbf{l}_r^{\delta}=n.$
\item[c)]  $\mathbf{l}_r^{\delta}\leq \min\{k_{\delta}-\sum_{i<r}
\mathbf{l}_i^{\delta}, \frac 1r \cdot (n-\sum_{i<r} i\cdot
\mathbf{l}_i^{\delta})\}$, for all $r\in[n]$.
\item[d)] If $k_{\delta}=1$, then $\mathbf{l}_n^{\delta}=1$ and
$\mathbf{l}_r^{\delta}=0$, for all $r\in [n-1]$.
\item[e)] If $k_{\delta}=n$, then $\mathbf{l}_1^{\delta}=n$ and
$\mathbf{l}_r^{\delta}=0$, for all $r\in [n] \setminus \{1\}$.
Specifically, $\delta=\epsilon$.
\end{enumerate}
\end{lmm}

{\bf Proof.} Assertions (a) and (b) are immediate from
definitions. Then, assertions (c), (d) and (e) are consequences of
the previous ones. \hfill $\Box$ \vspace{0.5cm}

Now, let us see that the number of Latin squares having a given
isotopism $\Theta\in\mathcal{I}_n$ in its autotopism group only
depends on the cycle structure of $\Theta$:

\begin{thr} \label{thr3} Let
$(\mathbf{l}_{\alpha},\mathbf{l}_{\beta},\mathbf{l}_{\gamma})$ be
the cycle structure of a Latin square isotopism and let us
consider
$\Theta_1=(\alpha_1,\beta_1,\gamma_1),\Theta_2=(\alpha_2,\beta_2,\gamma_2)\in
\mathcal{I}_n(\mathbf{l}_{\alpha},\mathbf{l}_{\beta},$
$\mathbf{l}_{\gamma})$. Then, $\Delta(\Theta_1)=\Delta(\Theta_2)$.
\end{thr}

{\bf Proof.}  Since $\Theta_1$ and $\Theta_2$ have the same cycle
structure, we can consider the isotopism
$\Theta=(\sigma_1,\sigma_2,\sigma_3)\in \mathcal{I}_n$, where:

\vspace{0.2cm}

\begin{enumerate}
\item[i)] $\sigma_1(c_{i,j}^{\alpha_1})=c_{i,j}^{\alpha_2}$, for all
$i\in \{0,1,...,k_{\alpha_1}\}$ and $j\in
\{0,1,...,\lambda_i^{\alpha_1}\}$,
\item[ii)] $\sigma_2(c_{i,j}^{\beta_1})=c_{i,j}^{\beta_2}$, for all
$i\in \{0,1,...,k_{\beta_1}\}$ and $j\in
\{0,1,...,\lambda_i^{\beta_1}\}$,
\item[iii)] $\sigma_3(c_{i,j}^{\gamma_1})=c_{i,j}^{\gamma_2}$, for all
$i\in \{0,1,...,k_{\gamma_1}\}$ and $j\in
\{0,1,...,\lambda_i^{\gamma_1}\}$.
\end{enumerate}

\vspace{0.2cm}

Now, let us see that $\Delta(\Theta_1)\leq \Delta(\Theta_2)$. If
$\Delta(\Theta_1)=0$, the result is immediate. Otherwise, let
$L_1=(l_{i,j})\in LS(\Theta_1)$ and let us see that
$L_1^{\Theta}=(l'_{i,j})\in LS(\Theta_2)$. Specifically, we must
prove that $(\alpha_2(i),\beta_2(j),\gamma_2(l'_{i,j}))\in
L_1^{\Theta}$, for all $(i,j,l'_{i,j})\in L_1^{\Theta}$. So, let
us consider $(i_0,j_0,l'_{i_0,j_0})\in L_1^{\Theta}$ and let
$r_0\in \{0,1,...,k_{\alpha_2}\}, u_0\in
\{0,1,...,\lambda_{r_0}^{\alpha_2}\}, s_0\in
\{0,1,...,k_{\beta_2}\}, v_0\in
\{0,1,...,\lambda_{s_0}^{\beta_2}\},$
$t_0\in\{0,1,...,k_{\gamma_2}\}$ and $w_0 \in
\{0,1,...,\lambda_{t_0}^{\gamma_2}\}$ be such that
$c_{r_0,u_0}^{\alpha_2}=i_0, c_{s_0,v_0}^{\beta_2}=j_0$ and
$c_{t_0,w_0}^{\gamma_2}=l'_{i_0,j_0}$. Thus:
$$(c_{r_0,u_0}^{\alpha_1}, c_{s_0,v_0}^{\beta_1},
c_{t_0,w_0}^{\gamma_1})=(\sigma_1^{-1}(i_0),\sigma_2^{-1}(j_0),\sigma_3^{-1}(l'_{i_0,j_0}))\in
L_1.$$ Next, since $L_1\in LS(\Theta)$, we have that:
$$(c_{r_0,u_0 + 1\ (mod\ \lambda_{r_0}^{\alpha_1})}^{\alpha_1},
c_{s_0,v_0+ 1\ (mod\ \lambda_{s_0}^{\beta_1})}^{\beta_1},
c_{t_0,w_0+ 1\ (mod\ \lambda_{t_0}^{\gamma_1})}^{\gamma_1})=$$
$$=(\alpha_1(c_{r_0,u_0}^{\alpha_1}),
\beta_1(c_{s_0,v_0}^{\beta_1}),
\gamma_1(c_{t_0,w_0}^{\gamma_1}))\in L_1.$$ Therefore,
$(\alpha_2(i_0),\beta_2(j_0),\gamma_2(l'_{i_0,j_0}))\in
L_1^{\Theta}$, because:
$$(c_{r_0,u_0 + 1\ (mod\ \lambda_{r_0}^{\alpha_2})}^{\alpha_2},
c_{s_0,v_0+ 1\ (mod\ \lambda_{s_0}^{\beta_2})}^{\beta_2},
c_{t_0,w_0+ 1\ (mod\ \lambda_{t_0}^{\gamma_2})}^{\gamma_2})=$$
$$=(\sigma_1(c_{r_0,u_0 + 1\ (mod\
\lambda_{r_0}^{\alpha_1})}^{\alpha_1}),\sigma_2(c_{s_0,v_0+ 1\
(mod\ \lambda_{s_0}^{\beta_1})}^{\beta_1}), \sigma_3(c_{t_0,w_0+
1\ (mod\ \lambda_{t_0}^{\gamma_1})}^{\gamma_1})).$$

Analogously, it is verified that
$L_2^{(\sigma_1^{-1},\sigma_2^{-1},\sigma_3^{-1})}\in
LS(\Theta_1)$, for all $L_2\in LS(\Theta_2)$, and hence, the
result follows. \hfill $\Box$ \vspace{0.5cm}

From Theorem \ref{thr3}, a catalogue of the cycle structures of
all possible autotopisms of a Latin square, which is the goal of
the present paper, seems to be useful, because it would simplify
the general calculus of the number $\Delta(\Theta)$, which is at
the moment an open problem. Now, in order to obtain the mentioned
catalogue, let us see some previous results.

\begin{prp} \label{prp1a} Let $\Theta=(\alpha,\beta,\gamma)\in \mathcal{I}_{n}$ be such that $\Delta(\Theta)>0$. If
$\mathbf{l}_n^{\alpha}=\mathbf{l}_n^{\beta}=\mathbf{l}_n^{\gamma}=1$,
then $n$ must be odd.
\end{prp}

{\bf Proof.} From Lemma \ref{lmm1}, $\alpha,\beta$ and $\gamma$
consist of a single $n$ cycle. Let
$\Theta'=(\alpha,\alpha,\alpha)\in \mathcal{I}_n$. The cycle
structure of $\Theta'$ is the same as that of $\Theta$ and,
therefore, from Theorem \ref{thr3},
$\Delta(\Theta')=\Delta(\Theta)>0$. Let $L\in LS(\Theta')$. By
definition, $L$ is a diagonally cyclic Latin square, which is
possible only if $n$ is odd (Theorem 6, {\bf \cite{Wanless04}}).
\hfill $\Box$ \vspace{0.5cm}

The following results are consequences of Theorems \ref{thr1} and
\ref{thr2}:

\begin{prp} \label{prp2} Let $\Theta\in \mathcal{I}_{n}(\mathbf{l}_{\alpha},\ \mathbf{l}_{\beta},\
\mathbf{l}_{\gamma})$ be such that $\Delta(\Theta)>0$. If there
exist $\delta\in\{\alpha,\beta,\gamma\}$ such that
$\mathbf{l}_1^{\delta}>0$, then it must be that
$\mathbf{l}_r^{\delta_1}=\mathbf{l}_r^{\delta_2}$, for all
$r\in[n]$, where $\delta_1$ and $\delta_2$ are the two
permutations in $\{\alpha,\beta,\gamma\}\setminus\{\delta\}$.
Specifically, if $\mathbf{l}_1^{\delta}>\lfloor \frac n2 \rfloor$,
then $\mathbf{l}_1^{\delta_1}=\mathbf{l}_1^{\delta_2}=0$.
\end{prp}

{\bf Proof.} For a given $\delta, \delta_1$ and $\delta_2$ in the
hypothesis, we will be in case (a) of Theorem \ref{thr1}, if
$\mathbf{l}_1^{\delta_1}>0$, or in case (b) of such a result, if
$\mathbf{l}_1^{\delta_1}=0$. In both cases, the two permutations
$\delta_1$ and $\delta_2$ must have the same cycle structure and,
therefore, it must be that
$\mathbf{l}_r^{\delta_1}=\mathbf{l}_r^{\delta_2}$, for all
$r\in[n]$. Specifically, if $\mathbf{l}_1^{\delta}>\lfloor \frac
n2 \rfloor$, we are in case (b) of Theorem \ref{thr1} and so, it
must be that $\mathbf{l}_1^{\delta_1}=\mathbf{l}_1^{\delta_2}=0$.
\hfill $\Box$ \vspace{0.5cm}

\begin{prp} \label{prp3} Let $n\geq 2$ and let $\Theta\in \mathcal{I}_{n}(\mathbf{l}_{\alpha},\ \mathbf{l}_{\beta},\
\mathbf{l}_{\gamma})$ be such that $\Delta(\Theta)>0$. If there
exist $\delta_1\in\{\alpha,\beta,\gamma\}$ and
$\delta_2\in\{\alpha,\beta,\gamma\}\setminus\{\delta_1\}$ such
that $\mathbf{l}_1^{\delta_1}\cdot \mathbf{l}_1^{\delta_2} >0$,
then the three permutations $\alpha$, $\beta$ and $\gamma$ have
the same cycle structure with at least one and at most $\lfloor
\frac n2 \rfloor$ fixed points. Specifically, it must be that
$1\leq
\mathbf{l}_1^{\alpha}=\mathbf{l}_1^{\beta}=\mathbf{l}_1^{\gamma}\leq
\lfloor \frac n2 \rfloor$ and $2\leq
k_{\alpha}=k_{\beta}=k_{\gamma}\leq \lfloor \frac n2 \rfloor +
\lfloor\frac {\lceil \frac n2 \rceil}2\rfloor$.
\end{prp}

{\bf Proof.} The first part of the lemma is immediate from Theorem
\ref{thr1}, because we would be in case (a) of that result.
Specifically, that theorem assures that $1\leq
\mathbf{l}_1^{\alpha}=\mathbf{l}_1^{\beta}=\mathbf{l}_1^{\gamma}\leq
\lfloor \frac n2 \rfloor$ and that
$k_{\alpha}=k_{\beta}=k_{\gamma}$. Now, since $\alpha,\beta$ and
$\gamma$ all have at least one fixed point, then they must have at
least two cycles, because $n\geq 2$. The upper bound of this
number of cycles is obtained when
$\mathbf{l}_1^{\alpha}=\mathbf{l}_1^{\beta}=\mathbf{l}_1^{\gamma}=
\lfloor \frac n2 \rfloor$ and the rest of the cycles have all of
them length $2$. \hfill $\Box$ \vspace{0.5cm}

\begin{prp} \label{prp4} Let $\Theta\in \mathcal{I}_{n}(\mathbf{l}_{\alpha},\ \mathbf{l}_{\beta},\
\mathbf{l}_{\gamma})$ be such that $\Delta(\Theta)>0$. If there
exists $t\in [n]$ such that $\mathbf{l}_t^{\gamma}>0$, then there
must exist $r,s\in[n]$ such that $\mathbf{l}_r^{\alpha}\cdot
\mathbf{l}_s^{\beta}>0$ and $t$ divides $l.c.m.(r,s)$.
\end{prp}

{\bf Proof.} Let $L=(l_{i,j})\in LS(\Theta)$ and let us consider
$t_0\in\{1,2,...,k_{\gamma}-1\}$ such that
$\lambda_{t_0}^{\gamma}=t$. Then, let
$r_0\in\{0,1,...,k_{\alpha}-1\}$, $s_0\in\{0,1,...,k_{\beta}-1\}$,
$u_0\in\{0,1...,\lambda_{r_0}^{\alpha}-1\}$ and
$v_0\in\{0,1...,\lambda_{s_0}^{\beta}-1\}$ be such that
$l_{c_{r_0,u_0}^{\alpha},c_{s_0,v_0}^{\beta}}=c_{t_0,0}^{\gamma}$.
Thus, from Theorem \ref{thr2}, $t=\lambda_{t_0}^{\gamma}$ must
divide $l.c.m.(\lambda_{r_0}^{\alpha},\lambda_{s_0}^{\beta})$.
Moreover, it is verified that
$\mathbf{l}_{\lambda_{r_0}^{\alpha}}^{\alpha}\geq 1 \leq
\mathbf{l}_{\lambda_{s_0}^{\beta}}^{\beta}$ and, therefore,
$\mathbf{l}_{\lambda_{r_0}^{\alpha}}^{\alpha}\cdot
\mathbf{l}_{\lambda_{s_0}^{\beta}}^{\beta}>0$. So, it is enough to
take $r=\lambda_{r_0}^{\alpha}$ and $s=\lambda_{s_0}^{\beta}$.
\hfill $\Box$ \vspace{0.3cm}

\begin{prp} \label{prp5} Let $\Theta\in \mathcal{I}_{n}(\mathbf{l}_{\alpha},\ \mathbf{l}_{\beta},\
\mathbf{l}_{\gamma})$ be such that $\Delta(\Theta)>0$. Let $r,s\in
[n]$ be such that $\mathbf{l}_r^{\alpha}\cdot
\mathbf{l}_s^{\beta}>0$ and let $m=l.c.m.(r,s)$. Then, there must
exist $t\in [m]$ such that:
\begin{enumerate}
\item[i)] $\mathbf{l}_t^{\gamma}>0$,
\item[ii)] $t$ divides $m$,
\item[iii)] $t$ does not divide any multiple of $r$ smaller than $m$,
\item[iv)] $t$ does not divide any multiple of $s$ smaller than
$m$.
\end{enumerate}
Indeed, if $g.c.d.(r,s)=1$, then it must be that $m\leq n$ and
$\mathbf{l}_m^{\gamma}>0$.
\end{prp}

{\bf Proof.} The result is an immediate consequence from Theorem
\ref{thr2}. \hfill $\Box$ \vspace{0.5cm}

Let $r,s\in [n]$ such that $\mathbf{l}_r^{\alpha}\cdot
\mathbf{l}_s^{\beta}>0$ and let us denote by $S_{r,s}^{\gamma}$
the set of $t$'s satisfying the four assertions of Proposition
\ref{prp5}. Finally, let us define the following sets:
$$S_{r,t}^{\beta}=\{u\in [n] : \mathbf{l}_u^{\beta}>0 \text{ and } S_{r,u}^{\gamma}=\{t\}\},$$
$$S_{s,t}^{\alpha}=\{u\in [n] : \mathbf{l}_u^{\alpha}>0 \text{ and
} S_{u,s}^{\gamma}=\{t\}\}.$$

Then, the following result is verified:

\begin{thr} \label{thr4} Let $\Theta\in \mathcal{I}_{n}(\mathbf{l}_{\alpha},\ \mathbf{l}_{\beta},\
\mathbf{l}_{\gamma})$ be such that $\Delta(\Theta)>0$. Let $t\in
[n]$ be such that $\mathbf{l}_t^{\gamma}>0$. Then, if $r,s\in [t]$
are such that $\mathbf{l}_r^{\alpha}>0$ and
$\mathbf{l}_s^{\beta}>0$, then it is verified that:
$$\sum_{u\in S_{r,t}^{\beta}} u\cdot \mathbf{l}_u^{\beta}\leq
t\cdot \mathbf{l}_t^{\gamma} \hspace{1cm} \text{ and }
\hspace{1cm} \sum_{u\in S_{s,t}^{\alpha}} u\cdot
\mathbf{l}_u^{\alpha}\leq t\cdot \mathbf{l}_t^{\gamma}.$$
\end{thr}

{\bf Proof.}  Let $L=(l_{i,j})\in LS(\Theta)$ and let us consider
$t_0\in\{0,1,...,k_{\gamma}-1\}$ such that
$\lambda_{t_0}^{\gamma}=t$. We will prove the result with the set
$S_{r,t}^{\beta}$, being analogous the proof with the set
$S_{s,t}^{\alpha}$. If $S_{r,t}^{\beta}=\emptyset$, then the
result is immediate. So, we can suppose that $S_{r,t}^{\beta}\neq
\emptyset$. Let $u\in S_{r,t}^{\beta}$ and let us consider
$r_0\in\{0,1,...,k_{\alpha}-1\}$ and
$u_0\in\{0,1,...,k_{\beta}-1\}$ such that
$\lambda_{r_0}^{\alpha}=r$ and $\lambda_{u_0}^{\beta}=u$. Since
$S_{r,u}^{\gamma}=\{t\}$, we have that, for all
$v\in\{0,1,...,u-1\}$, there must exist
$t_v\in\{0,1,...,k_{\gamma}-1\}$ such that
$\lambda_{t_v}^{\gamma}=t$ and $l_{c_{r_0,0}^{\alpha},
c_{u_0,v}^{\beta}}\in C_{t_v}^{\gamma}$. Therefore, as $L$ is a
Latin square, it must be that $u\cdot \mathbf{l}_u^{\beta}\leq
t\cdot \mathbf{l}_t^{\gamma}$. Since $u$ has been arbitrarily
taken in $S_{r,t}^{\beta}$, then, by working in the same
$(c_{r_0,0}^{\alpha} + 1)^{th}$ row of $L$, it must be that
$\sum_{u\in S_{r,t}^{\beta}} u\cdot \mathbf{l}_u^{\beta} \leq
t\cdot \mathbf{l}_t^{\gamma}$, because $L$ is a Latin square and
so, $L$ cannot have any repeated element in the mentioned row.
\hfill $\Box$ \vspace{0.5cm}

Let us see an example:

\begin{example}
\label{ex0} Let $\Theta\in \mathcal{I}_6((0,1,0,1,0,0),
(6,0,0,0,0,0), (0,1,0,1,0,0))$ and let us consider $r=t=4$. In
this case, $S_{4,4}^{\beta}=\{1\}$ and $\sum_{u\in
S_{4,4}^{\beta}} u\cdot \mathbf{l}_u^{\beta}=1\cdot 6 =
6 > 4 = 4 \cdot \mathbf{l}_4^{\gamma}$. Therefore, from Theorem
\ref{thr4}, it must be $\Delta(\Theta)=0$.
\end{example}

Let us observe that Theorem \ref{thr4} can be stated in a
conjugacy invariant way, by interchanging the role of rows,
columns and symbols. So, from Example \ref{ex0}, it can be deduced
that any isotopism with cycle structure $((0,1,0,1,0,0),
(0,1,0,1,0,0),(6,0,0,0,0,0))$ or $((6,0,0,0,0,0), (0,1,0,1,$
$0,0), (0,1,0,1,0,0))$ cannot be a Latin square autotopism.

\vspace{0.3cm}

Let us finish this section with a result corresponding to
autotopisms having cycles of prime lengths:

\begin{thr} \label{thr5} Let $\Theta\in \mathcal{I}_{n}(\mathbf{l}_{\alpha},\ \mathbf{l}_{\beta},\
\mathbf{l}_{\gamma})$ be such that $\mathbf{l}_p^{\alpha}\cdot
\mathbf{l}_p^{\beta}>0$, for some prime $p\in [n]$. If
$\mathbf{l}_1^{\gamma}<p\cdot
\max\{\mathbf{l}_p^{\alpha},\mathbf{l}_p^{\beta}\}$ and
$\mathbf{l}_p^{\gamma}=0$, then $\Delta(\Theta)=0$. Moreover, if
$\mathbf{l}_1^{\gamma}=0$ and
$\mathbf{l}_p^{\gamma}<\max\{\mathbf{l}_p^{\alpha},\mathbf{l}_p^{\beta}\}$,
then $\Delta(\Theta)=0$. Finally, if $p=2$,
$\mathbf{l}_1^{\gamma}=0$ and $\mathbf{l}_2^{\gamma}=1$, then
$\Delta(\Theta)=0$.
\end{thr}

{\bf Proof.} Let us suppose that $\Delta(\Theta)>0$ and let us
consider $L=(l_{i,j})\in LS(\Theta)$. We can suppose that
$\mathbf{l}_p^{\alpha}\leq \mathbf{l}_p^{\beta}$ (the reasoning is
similar in the other case). Let $p_0\in\{0,1,...,k_{\alpha}-1\}$
be such that $\lambda_{p_0}^{\alpha}=p$. Now, let us study each
part of the hypothesis:
\begin{enumerate}
\item[a)] Let us suppose that $\mathbf{l}_1^{\gamma}<p\cdot\max\{\mathbf{l}_p^{\alpha},\mathbf{l}_p^{\beta}\}=p\cdot \mathbf{l}_p^{\beta}$ and
$\mathbf{l}_p^{\gamma}=0$. From Theorem \ref{thr2}, since
$\mathbf{l}_p^{\gamma}=0$, we have that, for all
$p_1\in\{0,1,...,k_{\beta}-1\}$ such that
$\lambda_{p_1}^{\beta}=p$ and for all $v\in\{0,1,...,p-1\}$, it
must be that $l_{c_{p_0,0}^{\alpha},c_{p_1,v}^{\beta}}\in
Fix(\gamma)$. So, $\gamma$ must have at least $p\cdot
\mathbf{l}_p^{\beta}$ fixed points, because $L$ is a Latin square.
But then, we obtain a contradiction with being
$\mathbf{l}_1^{\gamma}<p\cdot\max\{\mathbf{l}_p^{\alpha},\mathbf{l}_p^{\beta}\}$.
So, it must be that $\Delta(\Theta)=0$.
\item[b)] Let us suppose that $\mathbf{l}_1^{\gamma}=0$ and
$\mathbf{l}_p^{\gamma}<\max\{\mathbf{l}_p^{\alpha},\mathbf{l}_p^{\beta}\}=\mathbf{l}_p^{\beta}$.
From Theorem \ref{thr2}, since $\mathbf{l}_1^{\gamma}=0$, we have
that, for all $p_1\in\{0,1,...,k_{\beta}-1\}$ such that
$\lambda_{p_1}^{\beta}=p$ and for all $v\in\{0,1,...,p-1\}$, there
must exist $t_{p_1,v}\in\{0,1,...,k_{\gamma}-1\}$ such that
$\lambda_{t_{p_1,v}}^{\gamma}=p$ and
$l_{c_{p_0,0}^{\alpha},c_{p_1,v}^{\beta}}\in
C_{t_{p_1,v}}^{\gamma}$. So, $\gamma$ must have at least $p\cdot
\mathbf{l}_p^{\beta}$ different elements in cycles of length $p$,
because $L$ is a Latin square. Specifically, $\gamma$ must have at
least $\mathbf{l}_p^{\beta}$ cycles of length $p$. But then, we
obtain a contradiction with being
$\mathbf{l}_p^{\gamma}<\max\{\mathbf{l}_p^{\alpha},\mathbf{l}_p^{\beta}\}$.
So, it must be that $\Delta(\Theta)=0$.
\item[c)] Let us suppose that $p=2$ and let us consider $\mathbf{l}_1^{\gamma}=0$ and
$\mathbf{l}_2^{\gamma}=1$. Let $p_1\in\{0,1,...,k_{\beta}-1\}$ be
such that $\lambda_{p_1}^{\beta}=2$ and let
$t\in\{0,1,...,k_{\gamma}-1\}$ be such that
$l_{c_{p_0,0}^{\alpha},c_{p_1,0}^{\beta}}\in C_t^{\gamma}$.  From
Theorem \ref{thr2}, $t$ must divide
$l.c.m.(\lambda_{p_0}^{\alpha},\lambda_{p_1}^{\beta})=2$. Then, it
must be that $t=2$, because $\mathbf{l}_1^{\gamma}=0$. Indeed, let
us observe that the four elements
$l_{c_{p_0,0}^{\alpha},c_{p_1,0}^{\beta}},
l_{c_{p_0,0}^{\alpha},c_{p_1,1}^{\beta}},
l_{c_{p_0,1}^{\alpha},c_{p_1,0}^{\beta}}$ and
$l_{c_{p_0,1}^{\alpha},c_{p_1,1}^{\beta}}$, must be in
$C_t^{\gamma}$, because $\mathbf{l}_2^{\gamma}=1$. Now, let
$w\in\{0,1\}$ be such that
$l_{c_{p_0,0}^{\alpha},c_{p_1,0}^{\beta}}=c_{t,w}^{\gamma}$. Then,
it must be that
$l_{c_{p_0,1}^{\alpha},c_{p_1,1}^{\beta}}=c_{t,w+1\ (mod\
2)}^{\gamma}$. Therefore, let us observe that
$l_{c_{p_0,0}^{\alpha},c_{p_1,1}^{\beta}}$ cannot be in
$C_t^{\gamma}$, because $L$ is a Latin square. So, we have a
contradiction and thus, it must be that $\Delta(\Theta)=0$. \hfill
$\Box$ \vspace{0.5cm}
\end{enumerate}

\section{\hspace{-0.4cm} Cycle structures of autotopisms of the Latin squares of order up to 11.}

All the results of the previous section have been implemented in a
computer program to generate all the possible cycle structures of
the set of non-trivial autotopisms of the Latin squares of order
up to $11$. We can see all these cycle structures in the below
tables. Let us observe that it is enough to show those autotopisms
$\Theta=(\alpha,\beta,\gamma)$ in which $k_{\alpha}\leq
k_{\beta}\leq k_{\gamma}$, because of the conjugacy of rows,
columns and symbols in Latin squares. Otherwise,
$(\mathbf{l}_{\alpha},\ \mathbf{l}_{\beta},\ \mathbf{l}_{\gamma})$
is a cycle structure of a Latin square autotopism if and only if
it can be found a permutation $\sigma\in S_3$ such that
$k_{\pi_{\sigma(0)}(\Theta)}\leq k_{\pi_{\sigma(1)}(\Theta)}\leq
k_{\pi_{\sigma(2)}(\Theta)}$ and
$(\mathbf{l}_{\pi_{\sigma(0)}(\Theta)},\
\mathbf{l}_{\pi_{\sigma(1)}(\Theta)},\
\mathbf{l}_{\pi_{\sigma(2)}(\Theta)})$ is a cycle structure of a
Latin square autotopism, where $\pi_i$ gives the $(i+1)^{th}$
component of $\Theta$, for all $i\in\{0,1,2\}$.

\vspace{0.25cm}

{\small
\begin{center} \begin{tabular}{|c|c|c|c|} \hline
$n$ & $\mathbf{l}_{\alpha}$ & $\mathbf{l}_{\beta}$ &
$\mathbf{l}_{\gamma}$ \\ \hline 2 & (0,1) & (0,1) & (2,0) \\
\hline 3 &  (0,0,1) & (0,0,1) & (0,0,1)\\
\cline{4-4} \ & \ & \ & (3,0,0) \\  \cline{2-4} \ & (1,1,0) &
(1,1,0) & (1,1,0)  \\ \hline 4 &  (0,0,0,1) & (0,0,0,1) &
(0,2,0,0) \\ \cline{4-4} \ & \ & \ & (2,1,0,0) \\ \cline{4-4} \ &
\ & \ &  (4,0,0,0) \\ \cline{2-4} \ & (0,2,0,0) & (0,2,0,0) &
(0,2,0,0) \\ \cline{4-4} \ & \ & \ & (2,1,0,0) \\ \cline{4-4} \ &
\ & \ & (4,0,0,0) \\ \cline{2-4} \ & (1,0,1,0) & (1,0,1,0) &
(1,0,1,0) \\ \cline{2-4} \ & (2,1,0,0) & (2,1,0,0) & (2,1,0,0) \\
\hline 5 & (0,0,0,0,1) & (0,0,0,0,1) & (0,0,0,0,1) \\ \cline{4-4}
\ & \ & \ & (5,0,0,0,0) \\ \cline{2-4} \ & (1,0,0,1,0) &
(1,0,0,1,0) & (1,0,0,1,0)) \\ \cline{2-4} \ & (1,2,0,0,0) &
(1,2,0,0,0) & (1,2,0,0,0) \\ \cline{2-4} \ & (2,0,1,0,0) &
(2,0,1,0,0) & (2,0,1,0,0) \\ \hline 6 & \ & \ & (0,0,2,0,0,0) \\
\cline{4-4} \ & \ & \ &
(1,1,1,0,0,0) \\
\cline{4-4} \ & \ & (0,0,0,0,0,1) & (2,2,0,0,0,0) \\
\cline{4-4} \ & (0,0,0,0,0,1) & \ & (3,0,1,0,0,0) \\
\cline{4-4} \ & \ & \ & (4,1,0,0,0,0) \\ \cline{4-4} \ & \ & \ &
(6,0,0,0,0,0) \\
\cline{3-4} \ & \ & (0,0,2,0,0,0) & (0,3,0,0,0,0) \\
\cline{2-4} \ & \ & \ & (0,0,2,0,0,0) \\
\cline{4-4} \ & (0,0,2,0,0,0) & (0,0,2,0,0,0) & (3,0,1,0,0,0) \\
\cline{4-4} \ & \ & \ & (6,0,0,0,0,0) \\
\cline{2-4} \ & (1,0,0,0,1,0) & (1,0,0,0,1,0) & (1,0,0,0,1,0) \\
\cline{2-4} \ & \ & \ & (2,2,0,0,0,0) \\
\cline{4-4} \ & (0,3,0,0,0,0) & (0,3,0,0,0,0) & (4,1,0,0,0,0) \\
\cline{4-4} \ & \ & \ & (6,0,0,0,0,0) \\ \cline{2-4} \ &
(2,0,0,1,0,0) &
(2,0,0,1,0,0) & (2,0,0,1,0,0) \\
\cline{2-4} \ & (2,2,0,0,0,0) & (2,2,0,0,0,0) & (2,2,0,0,0,0) \\
\cline{2-4} \ & (3,0,1,0,0,0) & (3,0,1,0,0,0) & (3,0,1,0,0,0) \\
\hline
\end{tabular}

\vspace{0.2cm}

Table 1: Cycle structures of non-trivial autotopisms of $LS(n)$,
for $2\leq n\leq 6$.
\end{center}}

\vspace{0.15cm}

{\small
\begin{center} \begin{tabular}{|c|c|c|c|} \hline
$n$ & $\mathbf{l}_{\alpha}$ & $\mathbf{l}_{\beta}$ & $\mathbf{l}_{\gamma}$ \\
\hline 7 & (0,0,0,0,0,0,1) & (0,0,0,0,0,0,1) & (0,0,0,0,0,0,1) \\
\cline{4-4} \ & \ & \ & (7,0,0,0,0,0,0) \\
\cline{2-4} \ & (1,0,0,0,0,1,0) & (1,0,0,0,0,1,0) &
(1,0,0,0,0,1,0) \\
\cline{2-4} \ & (1,0,2,0,0,0,0) & (1,0,2,0,0,0,0) &
(1,0,2,0,0,0,0) \\ \cline{2-4} \ & (1,1,0,1,0,0,0) &
(1,1,0,1,0,0,0) & (1,1,0,1,0,0,0) \\ \cline{2-4} \ &
(2,0,0,0,1,0,0) & (2,0,0,0,1,0,0) & (2,0,0,0,1,0,0) \\
\cline{2-4} \ & (1,3,0,0,0,0,0) & (1,3,0,0,0,0,0) &
(1,3,0,0,0,0,0) \\
\cline{2-4} \ & (3,0,0,1,0,0,0) & (3,0,0,1,0,0,0) &
(3,0,0,1,0,0,0) \\ \cline{2-4} \ & (3,2,0,0,0,0,0) &
(3,2,0,0,0,0,0) & (3,2,0,0,0,0,0)\\ \hline
\end{tabular}

\vspace{0.2cm}

Table 2: Cycle structures of non-trivial autotopisms of $LS(7)$.
\end{center}}

\vspace{0.3cm}

\begin{example} Let us consider $\Theta=((012345),(012)(345),(01)(23)(45))
\in \mathcal{I}_6((0,0,0,0,0,1), (0,0,2,$ $0,0,0),
(0,3,0,0,0,0))$. The following one is a Latin square of
$LS(\Theta)$: {\small $$\left(\begin{array}{cccccc}
  0 & 2 & 4 & 1 & 3 & 5 \\
  5 & 1 & 3 & 4 & 0 & 2 \\
  2 & 4 & 0 & 3 & 5 & 1 \\
  1 & 3 & 5 & 0 & 2 & 4 \\
  4 & 0 & 2 & 5 & 1 & 3 \\
  3 & 5 & 1 & 2 & 4 & 0 \\
\end{array}\right)$$}
\end{example}

\vspace{0.5cm}

\begin{example} Let us consider $\Theta=((01)(23)(45),(01)(23)(45),(01)(23)(45))$ $\in
\mathcal{I}_7((1,3,0,0,0,0,0),$
$(1,3,0,0,0,0,0),(1,3,0,0,0,0,0))$. The following one is a Latin
square of $LS(\Theta)$: {\small $$\left(\begin{array}{ccccccc}
  6 & 1 & 3 & 4 & 5 & 2 & 0\\
  0 & 6 & 5 & 2 & 3 & 4 & 1\\
  3 & 5 & 6 & 1 & 4 & 0 & 2\\
  4 & 2 & 0 & 6 & 1 & 5 & 3\\
  5 & 3 & 2 & 0 & 6 & 1 & 4\\
  2 & 4 & 1 & 3 & 0 & 6 & 5\\
  1 & 0 & 4 & 5 & 2 & 3 & 6\\
\end{array}\right)$$}
\end{example}

{\small
\begin{center} \begin{tabular}{|c|c|c|c|} \hline
$n$ & $\mathbf{l}_{\alpha}$ & $\mathbf{l}_{\beta}$ & $\mathbf{l}_{\gamma}$\\
\hline 8 & \ & \ & (0,0,0,2,0,0,0,0)\\
\cline{4-4} \ & \ & \ & (0,2,0,1,0,0,0,0)\\
\cline{4-4} \ & \ & \ & (0,4,0,0,0,0,0,0)\\
\cline{4-4} \ & \ & \ & (2,1,0,1,0,0,0,0)\\
\cline{4-4} \ & (0,0,0,0,0,0,0,1) & (0,0,0,0,0,0,0,1) & (2,3,0,0,0,0,0,0)\\
\cline{4-4} \ & \ & \ & (4,0,0,1,0,0,0,0)\\
\cline{4-4} \ & \ & \ & (4,2,0,0,0,0,0,0)\\
\cline{4-4} \ & \ & \ & (6,1,0,0,0,0,0,0)\\
\cline{4-4} \ & \ & \ & (8,0,0,0,0,0,0,0)\\
\cline{2-4} \ & \ & \ & (0,0,0,2,0,0,0,0)\\
\cline{4-4} \ & \ & \ & (0,2,0,1,0,0,0,0)\\
\cline{4-4} \ & \ & \ &
(0,4,0,0,0,0,0,0)\\
\cline{4-4} \ & \ & \ &
(2,1,0,1,0,0,0,0)\\
\cline{4-4} \ & (0,0,0,2,0,0,0,0) & (0,0,0,2,0,0,0,0) &
(2,3,0,0,0,0,0,0)\\
\cline{4-4} \ & \ & \ &
(4,0,0,1,0,0,0,0)\\
\cline{4-4} \ & \ & \ &
(4,2,0,0,0,0,0,0)\\
\cline{4-4} \ & \ & \ &
(6,1,0,0,0,0,0,0)\\
\cline{4-4} \ & \ & \ &
(8,0,0,0,0,0,0,0)\\
\cline{2-4} \ & (0,1,0,0,0,1,0,0) & (0,1,0,0,0,1,0,0) &
(2,0,0,0,0,1,0,0)\\
\cline{4-4} \ & \ & \ &
(2,0,2,0,0,0,0,0)\\
\cline{2-4} \ & (1,0,0,0,0,0,1,0) & (1,0,0,0,0,0,1,0) &
(1,0,0,0,0,0,1,0)\\
\cline{2-4} \ & \ & \ &
 (0,2,0,1,0,0,0,0)\\
\cline{4-4} \ & (0,2,0,1,0,0,0,0) & (0,2,0,1,0,0,0,0) &
 (2,1,0,1,0,0,0,0)\\
\cline{4-4} \ & \ & \ &
 (4,0,0,1,0,0,0,0)\\
\cline{2-4} \ & (2,0,0,0,0,1,0,0) & (2,0,0,0,0,1,0,0) &
 (2,0,0,0,0,1,0,0)\\
\cline{2-4} \ & \ & \ &
 (0,4,0,0,0,0,0,0)\\
\cline{4-4} \ & \ & \ &
 (2,3,0,0,0,0,0,0)\\
\cline{4-4} \ & (0,4,0,0,0,0,0,0) & (0,4,0,0,0,0,0,0) &
 (4,2,0,0,0,0,0,0)\\
\cline{4-4} \ & \ & \ &
 (6,1,0,0,0,0,0,0)\\
\cline{4-4} \ & \ & \ &
 (8,0,0,0,0,0,0,0)\\
\cline{2-4} \ & (2,0,2,0,0,0,0,0) & (2,0,2,0,0,0,0,0) &
 (2,0,2,0,0,0,0,0)\\
\cline{2-4} \ & (2,1,0,1,0,0,0,0) & (2,1,0,1,0,0,0,0) &
 (2,1,0,1,0,0,0,0))\\
\cline{2-4} \ & (3,0,0,0,1,0,0,0) & (3,0,0,0,1,0,0,0) &
 (3,0,0,0,1,0,0,0)\\
\cline{2-4} \ & (2,3,0,0,0,0,0,0) & (2,3,0,0,0,0,0,0) &
 (2,3,0,0,0,0,0,0)\\
\cline{2-4} \ & (4,0,0,1,0,0,0,0) & (4,0,0,1,0,0,0,0) &
 (4,0,0,1,0,0,0,0)\\
\cline{2-4} \ & (4,2,0,0,0,0,0,0) & (4,2,0,0,0,0,0,0) &
 (4,2,0,0,0,0,0,0)\\ \hline
\end{tabular}

\vspace{0.2cm}

Table 3: Cycle structures of non-trivial autotopisms of $LS(8)$.
\end{center}}

\vspace{0.2cm}

\begin{example} Let us consider $\Theta=((01)(23)(45)(67),(01)(23)(45)(67),(01)$ $(23)(45))
\in
\mathcal{I}_8((0,4,0,0,0,0,0,0),(0,4,0,0,0,0,0,0),(2,3,0,0,0,0,0,0))$.
The following one is a Latin square of $LS(\Theta)$: {\small
$$\left(\begin{array}{cccccccc}
  0 & 2 & 1 & 3 & 4 & 6 & 5 & 7\\
  3 & 1 & 2 & 0 & 6 & 5 & 7 & 4\\
  1 & 3 & 4 & 6 & 5 & 7 & 0 & 2\\
  2 & 0 & 6 & 5 & 7 & 4 & 3 & 1\\
  4 & 6 & 5 & 7 & 0 & 2 & 1 & 3\\
  6 & 5 & 7 & 4 & 3 & 1 & 2 & 0\\
  5 & 7 & 0 & 2 & 1 & 3 & 4 & 6\\
  7 & 4 & 3 & 1 & 2 & 0 & 6 & 5\\
\end{array}\right)$$}
\end{example}

\vspace{0.5cm}

{\small
\begin{center} \begin{tabular}{|c|c|c|c|} \hline
$n$ & $\mathbf{l}_{\alpha}$ & $\mathbf{l}_{\beta}$ & $\mathbf{l}_{\gamma}$\\
\hline 9 & \ & \ &
(0,0,0,0,0,0,0,0,1)\\
\cline{4-4} \ & \ & \ &
 (0,0,3,0,0,0,0,0,0)\\
\cline{4-4} \ & (0,0,0,0,0,0,0,0,1) & (0,0,0,0,0,0,0,0,1) &
 (3,0,2,0,0,0,0,0,0)\\
\cline{4-4} \ & \ & \ &
 (6,0,1,0,0,0,0,0,0)\\
\cline{4-4} \ & \ & \ &
 (9,0,0,0,0,0,0,0,0)\\
\cline{2-4} \ & \ & \ &
 (0,0,1,0,0,1,0,0,0)\\
\cline{4-4} \ & (0,0,1,0,0,1,0,0,0) & (0,0,1,0,0,1,0,0,0) &
 (0,3,1,0,0,0,0,0,0)\\
\cline{4-4} \ & \ & \ &
 (3,0,0,0,0,1,0,0,0)\\
\cline{4-4} \ & \ & \ &
 (3,3,0,0,0,0,0,0,0)\\
\cline{2-4} \ & (1,0,0,0,0,0,0,1,0) & (1,0,0,0,0,0,0,1,0) &
 (1,0,0,0,0,0,0,1,0)\\
\cline{2-4} \ & \ & \ &
 (0,0,3,0,0,0,0,0,0)\\
\cline{4-4} \ & (0,0,3,0,0,0,0,0,0) & (0,0,3,0,0,0,0,0,0) &
 (3,0,2,0,0,0,0,0,0)\\
\cline{4-4} \ & \ & \ &
 (6,0,1,0,0,0,0,0,0)\\
\cline{4-4} \ & \ & \ &
 (9,0,0,0,0,0,0,0,0)\\
\cline{2-4} \ & (1,0,0,2,0,0,0,0,0) & (1,0,0,2,0,0,0,0,0) &
 (1,0,0,2,0,0,0,0,0)\\
\cline{2-4} \ & (1,1,0,0,0,1,0,0,0) & (1,1,0,0,0,1,0,0,0) &
 (1,1,0,0,0,1,0,0,0)\\
\cline{2-4} \ & (2,0,0,0,0,0,1,0,0) & (2,0,0,0,0,0,1,0,0) &
 (2,0,0,0,0,0,1,0,0)\\
\cline{2-4} \ & (3,0,0,0,0,1,0,0,0) & (3,0,0,0,0,1,0,0,0) &
 (3,0,0,0,0,1,0,0,0)\\
\cline{2-4} \ & (1,4,0,0,0,0,0,0,0) & (1,4,0,0,0,0,0,0,0) &
 (1,4,0,0,0,0,0,0,0)\\
\cline{2-4} \ & (3,0,2,0,0,0,0,0,0) & (3,0,2,0,0,0,0,0,0) &
 (3,0,2,0,0,0,0,0,0)\\
\cline{2-4} \ & (4,0,0,0,1,0,0,0,0) & (4,0,0,0,1,0,0,0,0) &
 (4,0,0,0,1,0,0,0,0)\\
\cline{2-4} \ & (3,3,0,0,0,0,0,0,0) & (3,3,0,0,0,0,0,0,0) &
 (3,3,0,0,0,0,0,0,0)\\ \hline
\end{tabular}

\vspace{0.2cm}

Table 4: Cycle structures of non-trivial autotopisms of $LS(9)$.
\end{center}}

{\small
\begin{center} \begin{tabular}{|c|c|c|c|} \hline
$n$ & $\mathbf{l}_{\alpha}$ & $\mathbf{l}_{\beta}$ & $\mathbf{l}_{\gamma}$\\
\hline 10 & \ & \ & (0,0,0,0,2,0,0,0,0,0)\\
\cline{4-4} \ & \ & \ & (1,2,0,0,1,0,0,0,0,0)\\ \cline{4-4} \ & \
& \ & (3,1,0,0,1,0,0,0,0,0)\\ \cline{4-4} \ & \ & \ &
(2,4,0,0,0,0,0,0,0,0)\\ \cline{4-4} \ & (0,0,0,0,0,0,0,0,0,1) &
(0,0,0,0,0,0,0,0,0,1) & (5,0,0,0,1,0,0,0,0,0)\\ \cline{4-4} \ & \
& \ & (4,3,0,0,0,0,0,0,0,0)\\ \cline{4-4} \ & \ & \ &
(6,2,0,0,0,0,0,0,0,0)\\ \cline{4-4} \ & \ & \ &
(8,1,0,0,0,0,0,0,0,0)\\ \cline{4-4} \ & \ & \ &
(10,0,0,0,0,0,0,0,0,0)\\ \cline{3-4} \ & \ & (0,0,0,0,2,0,0,0,0,0)
& (0,5,0,0,0,0,0,0,0,0)\\ \cline{2-4} \ & \ & \ &
(0,0,0,0,2,0,0,0,0,0)\\ \cline{4-4} \ & (0,0,0,0,2,0,0,0,0,0) &
(0,0,0,0,2,0,0,0,0,0) & (5,0,0,0,1,0,0,0,0,0)\\ \cline{4-4} \ & \
& \ & (10,0,0,0,0,0,0,0,0,0)\\ \cline{2-4} \ &
(0,1,0,0,0,0,0,1,0,0) & (0,1,0,0,0,0,0,1,0,0) &
(2,0,0,0,0,0,0,1,0,0)\\ \cline{2-4} \ & (1,0,0,0,0,0,0,0,1,0) &
(1,0,0,0,0,0,0,0,1,0) & (1,0,0,0,0,0,0,0,1,0)\\ \cline{2-4} \ &
(0,1,0,2,0,0,0,0,0,0) & (0,1,0,2,0,0,0,0,0,0) &
(2,0,0,2,0,0,0,0,0,0)\\ \cline{2-4} \ & \ & \ &
(0,2,2,0,0,0,0,0,0,0)\\ \cline{4-4} \ & \ & \ &
(2,1,0,0,0,1,0,0,0,0)\\ \cline{4-4} \ & (0,2,0,0,0,1,0,0,0,0) &
(0,2,0,0,0,1,0,0,0,0) & (2,1,2,0,0,0,0,0,0,0)\\ \cline{4-4} \ & \
& \ & (4,0,0,0,0,1,0,0,0,0)\\ \cline{4-4} \ & \ & \ &
(4,0,2,0,0,0,0,0,0,0)\\ \cline{2-4} \ & (1,0,1,0,0,1,0,0,0,0) &
(1,0,1,0,0,1,0,0,0,0) & (1,0,1,0,0,1,0,0,0,0)\\ \cline{2-4} \ &
(2,0,0,0,0,0,0,1,0,0) & (2,0,0,0,0,0,0,1,0,0) &
(2,0,0,0,0,0,0,1,0,0)\\ \cline{2-4} \ & (1,0,3,0,0,0,0,0,0,0) &
(1,0,3,0,0,0,0,0,0,0) & (1,0,3,0,0,0,0,0,0,0)\\ \cline{2-4} \ &
(2,0,0,2,0,0,0,0,0,0) & (2,0,0,2,0,0,0,0,0,0) &
(2,0,0,2,0,0,0,0,0,0)\\ \cline{2-4} \ & (2,1,0,0,0,1,0,0,0,0) &
(2,1,0,0,0,1,0,0,0,0)  & (2,1,0,0,0,1,0,0,0,0) \\ \cline{2-4} \ &
(3,0,0,0,0,0,1,0,0,0) & (3,0,0,0,0,0,1,0,0,0)  &
(3,0,0,0,0,0,1,0,0,0) \\ \cline{2-4} \ & \ & \ &
(2,4,0,0,0,0,0,0,0,0) \\ \cline{4-4} \ & \ & \ &
(4,3,0,0,0,0,0,0,0,0) \\ \cline{4-4} \ & (0,5,0,0,0,0,0,0,0,0) &
(0,5,0,0,0,0,0,0,0,0)  & (6,2,0,0,0,0,0,0,0,0) \\ \cline{4-4} \ &
\ & \  & (8,1,0,0,0,0,0,0,0,0) \\ \cline{4-4} \ & \ & \  &
(10,0,0,0,0,0,0,0,0,0) \\ \cline{2-4} \ & (4,0,0,0,0,1,0,0,0,0) &
(4,0,0,0,0,1,0,0,0,0)  & (4,0,0,0,0,1,0,0,0,0) \\ \cline{2-4} \ &
(2,4,0,0,0,0,0,0,0,0) & (2,4,0,0,0,0,0,0,0,0)  &
(2,4,0,0,0,0,0,0,0,0) \\ \cline{2-4} \ & (4,0,2,0,0,0,0,0,0,0) &
(4,0,2,0,0,0,0,0,0,0)  & (4,0,2,0,0,0,0,0,0,0) \\ \cline{2-4} \ &
(5,0,0,0,1,0,0,0,0,0) & (5,0,0,0,1,0,0,0,0,0)  &
(5,0,0,0,1,0,0,0,0,0) \\ \cline{2-4} \ & (4,3,0,0,0,0,0,0,0,0) &
(4,3,0,0,0,0,0,0,0,0)  & (4,3,0,0,0,0,0,0,0,0) \\ \hline
\end{tabular}

\vspace{0.2cm}

Table 5: Cycle structures of non-trivial autotopisms of $LS(10)$.
\end{center}}

{\small
\begin{center} \begin{tabular}{|c|c|c|c|} \hline
$n$ & $\mathbf{l}_{\alpha}$ & $\mathbf{l}_{\beta}$ & $\mathbf{l}_{\gamma}$\\
\hline 11 & (0,0,0,0,0,0,0,0,0,0,1) & (0,0,0,0,0,0,0,0,0,0,1) & (0,0,0,0,0,0,0,0,0,0,1)\\
\cline{4-4} \ & \ & \ & (11,0,0,0,0,0,0,0,0,0,0)\\ \cline{2-4} \ &
(1,0,0,0,0,0,0,0,0,1,0) & (1,0,0,0,0,0,0,0,0,1,0) &
(1,0,0,0,0,0,0,0,0,1,0)\\ \cline{2-4} \ & (1,0,0,0,2,0,0,0,0,0,0)
& (1,0,0,0,2,0,0,0,0,0,0) & (1,0,0,0,2,0,0,0,0,0,0)\\ \cline{2-4}
\ & (1,1,0,0,0,0,0,1,0,0,0) & (1,1,0,0,0,0,0,1,0,0,0) &
(1,1,0,0,0,0,0,1,0,0,0)\\ \cline{2-4} \ & (2,0,0,0,0,0,0,0,1,0,0)
& (2,0,0,0,0,0,0,0,1,0,0) & (2,0,0,0,0,0,0,0,1,0,0)\\ \cline{2-4}
\ & (1,1,0,2,0,0,0,0,0,0,0) & (1,1,0,2,0,0,0,0,0,0,0) &
(1,1,0,2,0,0,0,0,0,0,0)\\ \cline{2-4} \ & (1,2,0,0,0,1,0,0,0,0,0)
& (1,2,0,0,0,1,0,0,0,0,0) & (1,2,0,0,0,1,0,0,0,0,0)\\ \cline{2-4}
\ & (2,0,1,0,0,1,0,0,0,0,0) & (2,0,1,0,0,1,0,0,0,0,0) &
(2,0,1,0,0,1,0,0,0,0,0)\\ \cline{2-4} \ & (3,0,0,0,0,0,0,1,0,0,0)
& (3,0,0,0,0,0,0,1,0,0,0) & (3,0,0,0,0,0,0,1,0,0,0)\\ \cline{2-4}
\ & (2,0,3,0,0,0,0,0,0,0,0) & (2,0,3,0,0,0,0,0,0,0,0) &
(2,0,3,0,0,0,0,0,0,0,0)\\ \cline{2-4} \ & (3,0,0,2,0,0,0,0,0,0,0)
& (3,0,0,2,0,0,0,0,0,0,0) & (3,0,0,2,0,0,0,0,0,0,0)\\ \cline{2-4}
\ & (4,0,0,0,0,0,1,0,0,0,0) & (4,0,0,0,0,0,1,0,0,0,0) &
(4,0,0,0,0,0,1,0,0,0,0)\\ \cline{2-4} \ & (1,5,0,0,0,0,0,0,0,0,0)
& (1,5,0,0,0,0,0,0,0,0,0) & (1,5,0,0,0,0,0,0,0,0,0)\\ \cline{2-4}
\ & (5,0,0,0,0,1,0,0,0,0,0) & (5,0,0,0,0,1,0,0,0,0,0) &
(5,0,0,0,0,1,0,0,0,0,0)\\ \cline{2-4} \ & (3,4,0,0,0,0,0,0,0,0,0)
& (3,4,0,0,0,0,0,0,0,0,0) & (3,4,0,0,0,0,0,0,0,0,0)\\ \cline{2-4}
\ & (5,0,2,0,0,0,0,0,0,0,0) & (5,0,2,0,0,0,0,0,0,0,0) &
(5,0,2,0,0,0,0,0,0,0,0)\\ \cline{2-4} \ & (5,3,0,0,0,0,0,0,0,0,0)
& (5,3,0,0,0,0,0,0,0,0,0) & (5,3,0,0,0,0,0,0,0,0,0)\\

\hline
\end{tabular}

\vspace{0.2cm}

Table 6: Cycle structures of non-trivial autotopisms of $LS(11)$.
\end{center}}

\vspace{0.5cm}

\begin{example} Let us consider $\Theta=((012345)(678), (012345)(678), (012)(34)$ $(56)(78))
\in \mathcal{I}_9((0,0,1,0,0,1,0,0,0),(0,0,1,0,0,1,0,0,0),
(0,3,1,0,0,0,0,$ $0,0))$. The following one is a Latin square of
$LS(\Theta)$:

{\small $$\left(\begin{array}{ccccccccc}
  0 & 2 & 1 & 4 & 6 & 8 & 3 & 5 & 7\\
  7 & 1 & 0 & 2 & 3 & 5 & 8 & 4 & 6\\
  6 & 8 & 2 & 1 & 0 & 4 & 5 & 7 & 3\\
  3 & 5 & 7 & 0 & 2 & 1 & 4 & 6 & 8\\
  2 & 4 & 6 & 8 & 1 & 0 & 7 & 3 & 5\\
  1 & 0 & 3 & 5 & 7 & 2 & 6 & 8 & 4\\
  4 & 7 & 5 & 3 & 8 & 6 & 0 & 2 & 1\\
  5 & 3 & 8 & 6 & 4 & 7 & 2 & 1 & 0\\
  8 & 6 & 4 & 7 & 5 & 3 & 1 & 0 & 2\\
\end{array}\right)$$}
\end{example}

\vspace{0.5cm}

\begin{example} Let us consider $\Theta=((012345)(678), (012345)(678), (012345)$ $(678))
\in \mathcal{I}_{10}((1,0,1,0,0,1,0,0,0,0),(1,0,1,0,0,1,0,0,0,0),
(1,0,1,0,0,1,$ $0,0,0,0))$. The following one is a Latin square of
$LS(\Theta)$:

{\small $$\left(\begin{array}{ccccccccccc}
  4 & 6 & 5 & 7 & 9 & 8 & 1 & 3 & 2 & 0\\
  6 & 5 & 7 & 0 & 8 & 9 & 3 & 2 & 4 & 1\\
  9 & 7 & 0 & 8 & 1 & 6 & 5 & 4 & 3 & 2\\
  7 & 9 & 8 & 1 & 6 & 2 & 4 & 0 & 5 & 3\\
  3 & 8 & 9 & 6 & 2 & 7 & 0 & 5 & 1 & 4\\
  8 & 4 & 6 & 9 & 7 & 3 & 2 & 1 & 0 & 5\\
  0 & 2 & 4 & 3 & 5 & 1 & 9 & 8 & 7 & 6\\
  2 & 1 & 3 & 5 & 4 & 0 & 8 & 9 & 6 & 7\\
  1 & 3 & 2 & 4 & 0 & 5 & 7 & 6 & 9 & 8\\
  5 & 0 & 1 & 2 & 3 & 4 & 6 & 7 & 8 & 9\\
\end{array}\right)$$}
\end{example}

\vspace{0.5cm}

\begin{example} Let us consider $\Theta=((01)(23)(45)(67), (01)(23)(45)(67),(01)$ $(23)(45)(67))
\in
\mathcal{I}_{11}((3,4,0,0,0,0,0,0,0,0,0),(3,4,0,0,0,0,0,0,0,0,0),(3,$
$4,0,0,$ $0,0,0,0,0,0,0))$. The following one is a Latin square of
$LS(\Theta)$:

{\small $$\left(\begin{array}{ccccccccccc}
  10 & 0 & 4 & 6 & 8 & 2 & 9 & 3 & 7 & 5 & 1\\
  1 & 10 & 7 & 5 & 3 & 8 & 2 & 9 & 6 & 4 & 0\\
  9 & 5 & 10 & 2 & 0 & 6 & 8 & 4 & 3 & 1 & 7\\
  4 & 9 & 3 & 10 & 7 & 1 & 5 & 8 & 2 & 0 & 6\\
  8 & 7 & 9 & 1 & 10 & 4 & 0 & 2 & 5 & 6 & 3\\
  6 & 8 & 0 & 9 & 5 & 10 & 3 & 1 & 4 & 7 & 2\\
  5 & 1 & 8 & 3 & 9 & 7 & 10 & 6 & 0 & 2 & 4\\
  0 & 4 & 2 & 8 & 6 & 9 & 7 & 10 & 1 & 3 & 5\\
  2 & 3 & 1 & 0 & 4 & 5 & 6 & 7 & 10 & 9 & 8\\
  7 & 6 & 5 & 4 & 2 & 3 & 1 & 0 & 8 & 10 & 9\\
  3 & 2 & 6 & 7 & 1 & 0 & 4 & 5 & 9 & 8 & 10\\
\end{array}\right)$$}
\end{example}

\section{Final remarks}

Apart from the previous cycles structures, the following ones
verify all the results of Section 3, although an exhaustive
computation proves that they do not correspond to any Latin square
autotopism:

{\small
\begin{center} $\begin{array}{|c|c|c|c|} \hline
n & \mathbf{l}_{\alpha} & \mathbf{l}_{\beta} & \mathbf{l}_{\gamma} \\
\hline 6 & (0, 0, 0, 0, 0, 1) & (0, 0, 0, 0, 0, 1) & (0, 3, 0, 0, 0, 0)\\
\cline{2-4} \ & (0, 1, 0, 1, 0, 0) & (0, 1, 0, 1, 0, 0) & (2, 0, 0, 1, 0, 0)\\
\cline{2-4} \ & (0, 3, 0, 0, 0, 0) & (0, 3, 0, 0, 0, 0) & (0, 3, 0, 0, 0, 0)\\
\hline 10 & (0, 0, 0, 0, 0, 0, 0, 0, 0, 1) & (0, 0, 0, 0, 0, 0,
0, 0, 0, 1) & (0, 5, 0, 0, 0, 0, 0, 0, 0, 0)\\
\cline{2-4} \ & (0, 2, 0, 0, 0, 1, 0, 0, 0, 0) & (0, 2, 0, 0,
0, 1, 0, 0, 0, 0) & (0, 2, 0, 0, 0, 1, 0, 0, 0, 0)\\
\cline{2-4} \ & (0, 5, 0, 0, 0, 0, 0, 0, 0, 0) & (0, 5, 0, 0, 0,
0, 0, 0, 0, 0) & (0, 5, 0, 0, 0, 0, 0, 0, 0, 0)\\\hline
\end{array}$ \end{center}}

Although in Section 4 we give all the cycle structures of
autotopisms of the Latin squares of order up to $11$, let us
remark that the properties of Section 3 can be implemented in an
algorithm to obtain all the cycle structures of autotopisms of the
Latin squares of greater orders.

\end{document}